\magnification 1200
\input plainenc
\input amssym
\input xepsf
\fontencoding{T2A}
\inputencoding{utf-8}
\tolerance 4000
\relpenalty 10000
\binoppenalty 10000
\parindent 1.5em

\hsize 17truecm
\vsize 24truecm
\hoffset 0truecm
\voffset -0.5truecm

\font\TITLE labx1440
\font\tenrm larm1000
\font\cmtenrm cmr10
\font\tenit lati1000
\font\tenbf labx1000
\font\teni cmmi10 \skewchar\teni '177
\font\tensy cmsy10 \skewchar\tensy '60
\font\tenex cmex10
\font\teneufm eufm10
\font\eightrm larm0800
\font\cmeightrm cmr8
\font\eightit lati0800
\font\eightbf labx0800
\font\eighti cmmi8 \skewchar\eighti '177
\font\eightsy cmsy8 \skewchar\eightsy '60
\font\eightex cmex8
\font\eighteufm eufm8

\font\cmsixrm cmr6

\font\sixbf labx0600
\font\sixi cmmi6 \skewchar\sixi '177
\font\sixsy cmsy6 \skewchar\sixsy '60
\font\sixeufm eufm6

\font\cmfiverm cmr5

\font\fivebf labx0500
\font\fivei cmmi5 \skewchar\fivei '177
\font\fivesy cmsy5 \skewchar\fivesy '60
\font\fiveeufm eufm5
\font\tencmmib cmmib10 \skewchar\tencmmib '177
\font\eightcmmib cmmib8 \skewchar\eightcmmib '177
\font\sevencmmib cmmib7 \skewchar\sevencmmib '177
\font\sixcmmib cmmib6 \skewchar\sixcmmib '177
\font\fivecmmib cmmib5 \skewchar\fivecmmib '177
\newfam\cmmibfam
\textfont\cmmibfam\tencmmib \scriptfont\cmmibfam\sevencmmib
\scriptscriptfont\cmmibfam\fivecmmib
\def\tenpoint{\def\rm{\fam0\tenrm}\def\it{\fam\itfam\tenit}%
	\def\bf{\fam\bffam\tenbf}
	\textfont0\cmtenrm \scriptfont0\cmsevenrm \scriptscriptfont0\cmfiverm
  	\textfont1\teni \scriptfont1\seveni \scriptscriptfont1\fivei
  	\textfont2\tensy \scriptfont2\sevensy \scriptscriptfont2\fivesy
  	\textfont3\tenex \scriptfont3\tenex \scriptscriptfont3\tenex
  	\textfont\itfam\tenit
	\textfont\bffam\tenbf \scriptfont\bffam\sevenbf
	\scriptscriptfont\bffam\fivebf
	\textfont\eufmfam\teneufm \scriptfont\eufmfam\seveneufm
	\scriptscriptfont\eufmfam\fiveeufm
	\textfont\cmmibfam\tencmmib \scriptfont\cmmibfam\sevencmmib
	\scriptscriptfont\cmmibfam\fivecmmib
	\normalbaselineskip 12pt
	\setbox\strutbox\hbox{\vrule height8.5pt depth3.5pt width0pt}%
	\normalbaselines\rm}
\def\eightpoint{\def\rm{\fam 0\eightrm}\def\it{\fam\itfam\eightit}%
	\def\bf{\fam\bffam\eightbf}%
	\textfont0\cmeightrm \scriptfont0\cmsixrm \scriptscriptfont0\cmfiverm
	\textfont1\eighti \scriptfont1\sixi \scriptscriptfont1\fivei
	\textfont2\eightsy \scriptfont2\sixsy \scriptscriptfont2\fivesy
	\textfont3\eightex \scriptfont3\eightex \scriptscriptfont3\eightex
	\textfont\itfam\eightit
	\textfont\bffam\eightbf \scriptfont\bffam\sixbf
	\scriptscriptfont\bffam\fivebf
	\textfont\eufmfam\eighteufm \scriptfont\eufmfam\sixeufm
	\scriptscriptfont\eufmfam\fiveeufm
	\textfont\cmmibfam\eightcmmib \scriptfont\cmmibfam\sixcmmib
	\scriptscriptfont\cmmibfam\fivecmmib
	\normalbaselineskip 11pt
	\abovedisplayskip 5pt
	\belowdisplayskip 5pt
	\setbox\strutbox\hbox{\vrule height7pt depth2pt width0pt}%
	\normalbaselines\rm
}

\def\No{\char 157}
\def\empty{}

\catcode`\@ 11
\catcode`\" 13
\def"#1{\ifx#1<\char 190\relax\else\ifx#1>\char 191\relax\else #1\fi\fi}

\def\newl@bel#1#2{\expandafter\def\csname l@#1\endcsname{#2}}
\openin 11\jobname .aux
\ifeof 11
	\closein 11\relax
\else
	\closein 11
	\input \jobname .aux
	\relax
\fi

\newcount\c@section
\newcount\c@subsection
\newcount\c@subsubsection
\newcount\c@equation
\newcount\c@bibl
\newcount\c@enum
\c@section 0
\c@subsection 0
\c@subsubsection 0
\c@equation 0
\c@bibl 0
\newdimen\d@enum
\d@enum=0pt
\def\lab@l{}
\def\label#1{\immediate\write 11{\string\newl@bel{#1}{\lab@l}}%
	\ifhmode\unskip\fi}
\def\eqlabel#1{\rlap{$(\equation)$}\label{#1}}

\def\section#1{\global\advance\c@section 1
	{\par\vskip 3ex plus 0.5ex minus 0.1ex
	\rightskip 0pt plus 1fill\leftskip 0pt plus 1fill\noindent
	{\bf\S\thinspace\number\c@section .~#1}\par\penalty 25000%
	\vskip 1ex plus 0.25ex}
	\gdef\lab@l{\number\c@section.}
	\c@subsection 0
	\c@subsubsection 0
	\c@equation 0
}
\def\subsection{\global\advance\c@subsection 1
	\par\vskip 1ex plus 0.1ex minus 0.05ex{\bf\number\c@subsection. }%
	\gdef\lab@l{\number\c@section.\number\c@subsection}%
	\c@subsubsection 0\c@equation 0%
}
\def\subsubsection{\global\advance\c@subsubsection 1
	\par\vskip 1ex plus 0.1ex minus 0.05ex%
	{\bf\number\c@subsection.\number\c@subsubsection. }%
	\gdef\lab@l{\number\c@section.\number\c@subsection.%
		\number\c@subsubsection}%
}
\def\equation{\global\advance\c@equation 1
	\gdef\lab@l{\number\c@section.\number\c@subsection.%
	\number\c@equation}{\rm\number\c@equation}
}
\def\bibitem#1{\global\advance\c@bibl 1
	[\number\c@bibl]%
	\gdef\lab@l{\number\c@bibl}\label{#1}
}
\def\ref@ref#1.#2:{\def\REF@{#2}\ifx\REF@\empty{\S\thinspace#1}%
	\else\ifnum #1=\c@section {#2}\else {\S\thinspace#1.#2}\fi\fi
}
\def\ref@eqref#1.#2.#3:{\ifnum #1=\c@section\ifnum #2=\c@subsection
	{(#3)}\else{#2\thinspace(#3)}\fi\else{\S\thinspace#1.#2\thinspace(#3)}\fi
}
\def\ref#1{\expandafter\ifx\csname l@#1\endcsname\relax
	{\bf ??}\else\edef\mur@{\csname l@#1\endcsname :}%
	{\expandafter\ref@ref\mur@}\fi
}
\def\eqref#1{\expandafter\ifx\csname l@#1\endcsname\relax
	{(\bf ??)}\else\edef\mur@{\csname l@#1\endcsname :}%
	{\expandafter\ref@eqref\mur@}\fi
}
\def\cite#1{\expandafter\ifx\csname l@#1\endcsname\relax
	{\bf ??}\else\hbox{\bf\csname l@#1\endcsname}\fi
}

\def\Im{\mathop{\rm Im}}
\def\im{\mathop{\rm im}}

\def\vraisup{\mathop{\rm vrai\;sup}}
\def\rank{\mathop{\rm rank}}
\def\Co{{\mathpalette\Wo@{}}C}
\def\Wo{{\mathpalette\Wo@{}}W}
\def\Wo@#1{\setbox0\hbox{$#1 W$}\dimen@\ht0\dimen@ii\wd0\raise0.65\dimen@%
\rlap{\kern0.35\dimen@ii$#1{}^\circ$}}

\catcode`\"=12
\def\bolddelta{\mathchar"0\hexnumber@\cmmibfam 0E}
\catcode`\"=13
\long\def\enumerate#1#2{%
	\setbox0\hbox{$#1^{\circ}.\ $}\d@enum=\wd0\global\advance\d@enum 2pt\c@enum=0%
	{\def\item{\global\advance\c@enum 1\par\hskip 0pt%
	\hbox to \d@enum{$\number\c@enum^{\circ}$.\hss}}%
	\par\smallskip #2\par\smallskip}
}
\def\proof{\par\medskip{\rm Д$\,$о$\,$к$\,$а$\,$з$\,$а$\,$т$\,$е$\,$л$\,$ь%
	$\,$с$\,$т$\,$в$\,$о.}\ }
\def\endproof{{\parfillskip 0pt\hfill$\square$\par}\medskip}

\catcode`\@ 12
\immediate\openout 11\jobname.aux


\frenchspacing\rm
\leftline{УДК~517.984}\vskip 5pt
{\TITLE\rightskip 0pt plus 1fill\leftskip 0pt plus 1fill\noindent
Об одном подходе к определению сингулярных дифференциальных операторов\par
\vskip 2ex\noindent\rm А.$\,$А.~Владимиров\footnote{}{\eightrm Работа поддержана
РФФИ, грант \No~16-01-00706.}\par}
\vskip 0.25cm
$$
	\vbox{\hsize 0.75\hsize\leftskip 0cm\rightskip 0cm
	\eightpoint\rm
	{\bf Аннотация:\/} На основе представления о тройках банаховых пространств
	даётся определение и характеризация основных свойств широкого класса
	граничных задач для обыкновенных дифференциальных уравнений произ\-вольного
	(в том числе нечётного) порядка с сингулярными коэффициентами.
	}
$$

\vskip 0.5cm

\section{Введение}
\subsection
В работе~\hbox{[\cite{Vl:2004}]} была указана конструкция, позволяющая дать
корректное определение ряда граничных задач для дифференциальных уравнений вида
$$
	\sum_{k=0}^n(-1)^{n-k}(p_ky^{(n-k)})^{(n-k)}=f
$$
с сингулярными коэффициентами $p_k\in W_2^{-k}[0,1]$. В случае $n=2$ при этом было
показано, что действие соответствующих неограниченных операторов в пространстве
$L_2[0,1]$ допускает описание посредством систем обычных дифференциальных уравнений
для абсолютно непрерывных функций. Применительно к задачам более высоких порядков
аналогичные представления явным образом не указывались, хотя возможность такого
указания в свете развитой теории была совершенно прозрачной. Значимость конструкции
из~[\cite{Vl:2004}] может быть продемонстрирована, например, развитым на её основе
в работе~[\cite{Vl:2016}] простым подходом к изучению осцилляционных свойств
собственных функций сингулярных обыкновенных дифференциальных операторов~---
включая ряд так называемых "<многоточечных"> задач, рассмотрение которых ранее
обычно проводилось с использованием весьма трудоёмких косвенных методов
(см.~имеющиеся в~[\cite{Vl:2016}] ссылки).

Идеологические основы подхода работы~[\cite{Vl:2004}] были заложены в более ранней
работе~[\cite{NSh:1999}], где схожим способом изучались дифференциальные уравнения
второго порядка.

В последнее время интерес к этой тематике возобновился.
Так, в работе~[\cite{MSh:2016}] были~--- вне связи с подходом
работы~[\cite{Vl:2004}]~--- предложены регуляризованные представления для формально
несколько более широкого, нежели рассмотренный в~[\cite{Vl:2004}], класса
дифференциальных уравнений чётного порядка. Поэтому представляется не лишённым
интереса вопрос о расширении границ применения использованной в~[\cite{Vl:2004}]
методики. Характеризации некоторых возможных направлений такого расширения
и посвящается настоящая статья.

\subsection
Структура статьи имеет следующий вид. В~\ref{par:Sob} даётся определение
и указываются важнейшие свойства вспомогательных функциональных пространств,
используемых далее при описании основных объектов изучения. В~\ref{par:Oper}
непосредственно определяются и исследуются операторы, отвечающие допустимым
в рамках развиваемой теории граничным задачам для дифференциальных уравнений
с сингулярными коэффициен\-тами. Наконец, в~\ref{par:Ekz} приводятся примеры
приложения развитой теории к ряду конкретных ситуаций.

Все рассматриваемые далее линейные пространства предполагаются комплексными.

Наименьшим натуральным числом на всём протяжении статьи считается $0$.


\section{Квазидифференциальные пространства Соболева}\label{par:Sob}
\subsection
Пусть $A$~--- такая система функций класса $L_1[0,1]$, что для каждой пары
натуральных индексов $i$ и $j\leqslant i+1$ найдётся соответствующая ей функция
$A_{ij}\in L_1[0,1]$. Символом $C_A^n[0,1]$ мы далее будем обозначать
подпространство в пространстве $\{C[0,1]\}^{n+1}$ непрерывных вектор-функций
с $n+1$~компонентами, выделенное системой понимаемых в смысле обобщённого
дифференцирования уравнений
$$
	Y_i'=\sum_{j=0}^{i+1}A_{ij}Y_j,\qquad i<n.\leqno(\equation)
$$\label{eq:sob}%
Символом $W_{s,A}^n[0,1]$, где $n>0$ и $s\in [1,\infty)$, мы будем обозначать
результат пополнения прост\-ранства $C_A^n[0,1]$ по норме
$$
	\|Y\|_{W_{s,A}^n[0,1]}\rightleftharpoons\sum_{i=0}^{n-1}\|Y_i\|_{C[0,1]}+
		\left(\int_0^1|A_{n-1,n}|\cdot|Y_n|^s\,dx\right)^{1/s}.
$$
Наконец, символами $\Co_A^n[0,1]$ и $\Wo_{s,A}^n[0,1]$ мы будем обозначать
подпространства, получаемые замыканием в пространствах $C_A^n[0,1]$
и $W_{s,A}^n[0,1]$, соответственно, множеств принадлежащих этим пространствам
финитных на интервале $(0,1)$ вектор-функций. Тривиальным образом имеет место
следующий факт:

\subsubsection\label{prop:1.1.1}
{\it Если при каждом $i<n$ первообразная функции $|A_{i,i+1}|$ строго монотонна,
то естественное вложение $Y\mapsto Y_0$ пространства $W_{1,A}^n[0,1]$
в пространство $C[0,1]$ является инъективным.
}

\medskip
Обычные соболевские пространства $W_s^n[0,1]$ и $\Wo_s^n[0,1]$, где $s\in
[1,\infty)$, в рамках изложенной конструкции отвечают (с точностью до выбора
эквивалентной нормы) ситуации
$$
	A_{ij}(x)\equiv\cases{1&при $j=i+1$,\cr 0&иначе.}
$$

\subsection
Обозначим символом $M_n$ решение понимаемой в смысле обобщённого
дифферен\-цирования начальной задачи
$$
	(M_n)_{ij}'=\kern -1em\sum_{k=0}^{\inf\{n-1,i+1\}}\kern -1em
		A_{ik}\cdot (M_n)_{kj},\qquad
		(M_n)_{ij}(0)=\cases{1&при $i=j$,\cr 0&иначе,}
$$
где каждый из индексов $i$ и $j$ пробегает промежуток $\{0,\ldots,n-1\}$.
Очевидным образом, каждая из функций $(M_n)_{ij}$ принадлежит классу $W_1^1[0,1]$.
Сопоставляя в каждой точке $x\in [0,1]$ имеющей размер $n\times n$ матрице
элементов $(M_n)_{ij}(x)$ обратную матрицу элементов $(M_n^{-1})_{ij}(x)$,
получаем набор также принадлежащих классу $W_1^1[0,1]$ функций $(M_n^{-1})_{ij}$.
При этом всякая вектор-функция $Y\in W_{1,A}^n[0,1]$ восстанавливается по своей
последней компоненте согласно правилам
$$
	Y_i(x)=\sum_{j=0}^{n-1}(M_n)_{ij}(x)\cdot\left[Y_j(0)+
		\int_0^x A_{n-1,n}(t)\cdot (M_n^{-1})_{j,n-1}(t)Y_n(t)\,dt\right],
		\qquad i<n.\leqno(\equation)
$$\label{eq:Ysob}%

\subsubsection\label{prop:1.2.1}
{\it Пусть при некотором $n>0$ первообразная функции $|A_{n-1,n}|$ строго
монотонна. Тогда естественное вложение пространства $C_A^n[0,1]$ в пространство
$C_A^{n-1}[0,1]$ имеет плотный образ.
}%
\proof
Зафиксируем произвольную вектор-функцию $Y\in C_A^{n-1}[0,1]$ и поста\-вим ей
в соответствие вектор-функцию $\Phi\in\{C[0,1]\}^n$ вида
$$
	\Phi_i(x)\equiv Y_i(x)-\int_0^x\left[\sum_{j=0}^{\inf\{n-1,i+1\}}
		\kern -1em A_{ij}(t)Y_j(t)\right]\,dt,\qquad i<n.
$$
Зафиксируем также последовательность $\{Y_\alpha\}_{\alpha=0}^\infty$
вектор-функций класса $C_A^n[0,1]$, удовлетво\-ряющих при $i<n$ равенствам
$(Y_\alpha)_i(0)=Y_i(0)$ и таких, что первообразные суммируемых функций
$A_{n-1,n}\cdot (Y_\alpha)_n$ равномерно стремятся при $\alpha\to\infty$ к функции
$\Phi_{n-1}$. Из пред\-ставления~\eqref{eq:Ysob} и соотношений $(M_n^{-1})_{j,n-1}
\in W_1^1[0,1]$ тогда немедленно вытекает факт равномерной сходимости
функциональных последовательностей $\{(Y_\alpha)_i\}_{\alpha=0}^\infty$,
где $i<n$. Тривиальные тождества
$$
	Y_{n-1}(0)+\int_0^x A_{n-1,n}(t)(Y_\alpha)_n(t)\,dt=(Y_\alpha)_{n-1}(x)-
		\int_0^x\left[\sum_{j=0}^{n-1} A_{n-1,j}(t)(Y_\alpha)_j(t)
		\right]\,dt
$$
вместе с известными общими свойствами интегральных уравнений первого рода
озна\-чают теперь, что пределом последовательности образов вектор-функций
$Y_\alpha\in C_A^n[0,1]$ является в точности исходная вектор-функция
$Y\in C_A^{n-1}[0,1]$.
\endproof

\subsubsection\label{prop:mallev}
{\it Пусть при некотором $n>0$ первообразная функции $|A_{n-1,n}|$ строго
монотонна, а набор суммируемых функций $\{f_j\}_{j=0}^n$ удовлетворяет тождеству
$$
	(\forall Y\in\Co_A^n[0,1])\qquad
		\sum_{j=0}^n\int_0^1 f_j\overline{Y_j}\,dx=0.\leqno(\equation)
$$\label{eq:00ant}%
Тогда существует функция $h\in W_1^1[0,1]$, при почти всех $x\in [0,1]$
удовлетворяющая равен\-ствам $f_n(x)=\overline{A_{n-1,n}(x)}\,h(x)$. При этом
справедливо также тождество
$$
	(\forall Y\in\Co_A^{n-1}[0,1])\qquad
		\sum_{j=0}^{n-1}\int_0^1 g_j\overline{Y_j}\,dx=0,
$$
где положено
$$
	g_j\rightleftharpoons\cases{f_j-\overline{A_{n-1,j}}\,h,&при $j<n-1$,\cr
		f_{n-1}-\overline{A_{n-1,n-1}}\,h-h',&при $j=n-1$.}
$$
}%
\proof
С учётом представления~\eqref{eq:Ysob}, заведомо найдётся функция
$\varphi\in W_1^1[0,1]$, удовлетворяющая тождеству
$$
	(\forall Y\in\Co_A^n[0,1])\qquad\int_0^1 [f_n+
		\overline{A_{n-1,n}}\,\varphi]\,\overline{Y_n}\,dx=0.
$$
Функция $f_n+\overline{A_{n-1,n}}\varphi\in L_1[0,1]$ при этом заведомо принадлежит
линейной оболочке набора $\{\overline{A_{n-1,n}\cdot
(M_n^{-1})_{j,n-1}}\}_{j=0}^{n-1}$, что автоматически означает существование
функции $h\in W_1^1[0,1]$ с требуемым свойством. Выражая
теперь, с использованием определе\-ния~\eqref{eq:sob}, в тождестве~\eqref{eq:00ant}
суммируемую функцию $A_{n-1,n}Y_n$ через функции набора $\{Y_i\}_{i=0}^{n-1}$
и интегрируя по частям, убеждаемся в справедливости соотношения
$$
	(\forall Y\in\Co_A^n[0,1])\qquad
		\sum_{j=0}^{n-1}\int_0^1 g_j\overline{Y_j}\,dx=0.
$$
Учёт утверждения~\ref{prop:1.2.1} завершает доказательство.
\endproof

\subsection
Символом $Y^\wedge\in\Bbb C^{2n}$, где $Y\in W_{1,A}^n[0,1]$, мы далее будем
обозначать вектор гранич\-ных значений
$$
	Y^\wedge_k\rightleftharpoons\cases{Y_k(0)&при $k<n$,\cr Y_{k-n}(1)&иначе.}
$$
На основе произвольно фиксированной матрицы $U\in\Bbb C^{2n\times 2n}$
могут быть определены подпространства $C_{A,U}^n[0,1]$ и $W_{s,A,U}^n[0,1]$,
выделенные, соответственно, внутри $C_A^n[0,1]$ и $W_{s,A}^n[0,1]$ системой
граничных условий $UY^\wedge=0$.

На указанных подпространствах могут быть заданы полулинейные функционалы вида
$$
	\langle F,Y\rangle\equiv\cases{\displaystyle\sum_{i=0}^n
		\int_0^1 f_i\overline{Y_i}\,dx&при $Y\in C_{A,U}^n[0,1]$,\cr
		\displaystyle\sum_{i=0}^{n-1}\int_0^1 f_i\overline{Y_i}\,dx+
		\int_0^1 |A_{n-1,n}|^{1/s}\,f_n\overline{Y_n}\,dx&
		при $Y\in W_{s,A,U}^n[0,1]$,}\leqno(\equation)
$$\label{eq:funkc}%
где $f_i\in L_1[0,1]$. В случае рассмотрения пространства $W_{s,A,U}^n[0,1]$ здесь
дополнительно предполагается, что выполнено условие $f_n\in L_{s/(s-1)}[0,1]$,
причём для почти любого $x\in [0,1]$ из равенства $A_{n-1,n}(x)=0$ следует
равенство $f_n(x)=0$. Первообразные всех функций $|A_{i,i+1}|$, где $i<n$,
мы на протяжении настоящего пункта считаем строго монотонными.
Из представления~\eqref{eq:Ysob} немедленно вытекает, что всякий функционал
вида~\eqref{eq:funkc} может быть переписан в форме
$$
	\langle F,Y\rangle\equiv\cases{\displaystyle\int_0^1 g\overline{Y_n}\,dx+
		\sum_{i=0}^{n-1}\mu_i\overline{Y_i(0)}&при $Y\in C_{A,U}^n[0,1]$,\cr
		\displaystyle\int_0^1 |A_{n-1,n}|^{1/s}\,g\overline{Y_n}\,dx+
		\sum_{i=0}^{n-1}\mu_i\overline{Y_i(0)}&при $Y\in W_{s,A,U}^n[0,1]$,}
	\leqno(\equation)
$$\label{eq:funkc1}%
где $g\in L_1[0,1]$ в случае пространства $C_{A,U}^n[0,1]$,
и $g\in L_{s/(s-1)}[0,1]$ в случае пространства $W_{s,A,U}^n[0,1]$. Из того же
представления легко получается, что всякий функционал вида~\eqref{eq:funkc1}
допускает обратную перезапись в форме~\eqref{eq:funkc}. Сказанное означает
замкнутость линейного множества функционалов вида~\eqref{eq:funkc} в сопряжённом
к рассматриваемому пространстве. Из известной теоремы об общем виде полулинейного
непрерывного функционала на лебеговском пространстве легко выводится также,
что в случае пространства $W_{s,A,U}^n[0,1]$ множество функционалов
вида~\eqref{eq:funkc} в точности совпадает с сопряжённым пространством.

Пространство заданных на $C_{A,U}^n[0,1]$ функционалов вида~\eqref{eq:funkc}
мы будем далее обоз\-начать символом $W_{1,A,U}^{-n}[0,1]$. Аналогичным образом,
пространство таких функциона\-лов, заданных на $W_{s,A,U}^n[0,1]$, мы будем
обозначать символом $W_{s/(s-1),A,U}^{-n}[0,1]$.


\section{Сингулярные дифференциальные операторы}\label{par:Oper}
\subsection
Пусть $B$~--- система функций, обладающая аналогичными системе $A$ свойствами
и такая, что при почти всех $x\in [0,1]$ равенства $A_{n-1,n}(x)=0$
и $B_{m-1,m}(x)=0$ равносильны. Пусть также зафиксирован параметр
$s\in [1,+\infty)$, две матрицы $V\in C^{2m\times 2m}$,
$Q\in\Bbb C^{2m\times 2n}$ и система функций со следующими свойствами:

\medskip
\item{$\bullet$} $p_{nm},p_{nm}^{-1}\in L_\infty[0,1]$.

\smallskip
\item{$\bullet$} $p_{im}\in L_{s/(s-1)}[0,1]$ при $i<n$, причём в случае $s\neq 1$
для почти любого $x\in [0,1]$ из равенства $B_{m-1,m}(x)=0$ следует равенство
$p_{im}(x)=0$.

\smallskip
\item{$\bullet$} $p_{nj}\in L_s[0,1]$ при $j<m$, причём для почти любого
$x\in [0,1]$ из равенства $A_{n-1,n}(x)=0$ следует равенство $p_{nj}(x)=0$.

\smallskip
\item{$\bullet$} $p_{ij}\in L_1[0,1]$ при $i<n$ и $j<m$.

\medskip\noindent
Тогда может быть задан ограниченный оператор $T\colon W_{s,A,U}^n[0,1]\to
W_{s,B,V}^{-m}[0,1]$ вида
$$
	\displaylines{\hbox to \displaywidth{\eqlabel{eq:kvdif}\kern 1truecm
		$\displaystyle\langle TY,Z\rangle\equiv\int_0^1 |A_{n-1,n}|^{1/s}
		\cdot |B_{m-1,m}|^{(s-1)/s}\,p_{nm}Y_n\overline{Z_m}\,dx+
		{}$\hfill}\cr
		{}+\sum_{i=0}^{n-1}\int_0^1 |B_{m-1,m}|^{(s-1)/s}\,p_{im}Y_i
		\overline{Z_m}\,dx+\sum_{j=0}^{m-1}\int_0^1 |A_{n-1,n}|^{1/s}\,
		p_{nj}Y_n\overline{Z_j}\,dx+{}\cr
		\kern 4truecm{}+\sum_{i=0}^{n-1}\sum_{j=0}^{m-1}\int_0^1 p_{ij}Y_i
		\overline{Z_j}\,dx+\langle QY^\wedge,Z^\wedge\rangle.}
$$

\subsubsection
{\it Всякий оператор $T$ вида~\eqref{eq:kvdif} фредгольмов индекса $n-m-\rank U+
\rank V$.
}%
\proof
Справедливость доказываемого утверждения легко устана\-вливается на основе
представления~\eqref{eq:Ysob}, а также факта наличия у оператора
$\hat T\colon L_s([0,1];\,|A_{n-1,n}|)\to L_s([0,1];\,|B_{m-1,m}|)$ вида 
$$
	\hat Ty\rightleftharpoons\left|{A_{n-1,n}\over B_{m-1,m}}\right|^{1/s}\,
		p_{nm}y
$$
обратного с оценкой нормы $\vraisup_{x\in [0,1]}|p_{nm}^{-1}(x)|$.
\endproof

\subsection
Результаты предыдущего параграфа показывают, что в случае строгой монотонно\-сти
первообразных функций $|B_{j,j+1}|$ уравнения $TY=F$ допускают эквивалентную запись
в форме граничных задач для систем дифференциальных уравнений для абсо\-лютно
непрерывных функций. А именно, для искомого решения $Y\in W_{s,A,U}^n[0,1]$
с очевидностью могут быть определены квазипроизводные $y^{[i]}\rightleftharpoons
Y_i\in W_1^1[0,1]$, где $i<n$, удовлетво\-ряющие системе уравнений
$$
	(y^{[i]})'=\sum_{j=0}^{i+1}A_{ij}y^{[j]},\qquad i<n-1.\leqno(\equation)
$$\label{eq:kvdif1}%
Далее, из утверждения~\ref{prop:mallev} вытекает существование квазипроизводной
$y^{[n]}\in W_1^1[0,1]$, удовлетворяющей соотношению
$$
	|A_{n-1,n}|^{1/s}Y_n=p_{nm}^{-1}\cdot\left[{\overline{B_{m-1,m}}\over
		|B_{m-1,m}|^{(s-1)/s}}\cdot y^{[n]}-\sum_{i=0}^{n-1}p_{im}y^{[i]}
		+f_m\right],\leqno(\equation)
$$\label{eq:kvdifYn}%
а тогда и уравнению
$$
	\leqalignno{(y^{[n-1]})'&=\sum_{i=0}^{n-1}\left[A_{n-1,i}-{p_{nm}^{-1}p_{im}\,
		A_{n-1,n}\over|A_{n-1,n}|^{1/s}}\right]\cdot y^{[i]}+{}&
		\eqlabel{eq:kvdif2}\cr &\kern 3truecm{}+
		{p_{nm}^{-1}A_{n-1,n}\overline{B_{m-1,m}}\over |A_{n-1,n}|^{1/s}\,
		|B_{m-1,m}|^{(s-1)/s}}\cdot y^{[n]}+{p_{nm}^{-1}A_{n-1,n}f_m\over
		|A_{n-1,n}|^{1/s}}.}
$$
Наконец, утверждение~\ref{prop:mallev} вместе с соотношением~\eqref{eq:kvdifYn}
означают существование квазипроизводных $y^{[n+m-j-1]}\in W_1^1[0,1]$, где $j<m-1$,
подчиняющихся уравнениям
$$
	\displaylines{\hbox to \displaywidth{\eqlabel{eq:kvdif3}\kern 1truecm
		$\displaystyle (y^{[n+m-j-1]})'=\sum_{i=0}^{n-1}\left[p_{ij}-
		p_{nm}^{-1}p_{im}p_{nj}\right]\cdot y^{[i]}+{}$\hfill}\cr
		{}+\left[p_{nm}^{-1}
		\cdot{\overline{B_{m-1,m}}\over |B_{m-1,m}|^{(s-1)/s}}\cdot p_{nj}-
		\overline{B_{m-1,j}}\right]\cdot y^{[n]}-{}\cr
		\kern 4truecm{}-\sum_{k=\sup\{j,1\}}^{m-1}\overline{B_{k-1,j}}
		y^{[n+m-k]}-f_j+p_{nm}^{-1}p_{nj}f_m,\qquad j<m.}
$$
При этом, очевидным образом, оказывается справедливым тождество
$$
	\langle TY-F,Z\rangle\equiv\langle QY^\wedge-Y^\vee,Z^\wedge\rangle,
$$
где положено
$$
	Y^\vee_k\rightleftharpoons\cases{y^{[n+m-k-1]}(0)&при $k<m$,\cr
		-y^{[n+2m-k-1]}(1)&иначе.}
$$
Соответственно, исходное уравнение $TY=F$ равносильно системе дифференциальных
урав\-нений~\eqref{eq:kvdif1}, \eqref{eq:kvdif2} и~\eqref{eq:kvdif3}, рассмотренной
совместно с набором граничных условий
$$
	UY^\wedge=0,\qquad QY^\wedge-Y^\vee\in\im V^*.
$$

\subsection
Пусть теперь для рассматриваемых пространств $W_{s,A,U}^n[0,1]$
и $W_{s,B,V}^{-m}[0,1]$, а также некоторого гильбертова пространства $\frak H$,
зафиксированы вложения $I\colon W_{s,A,U}^n[0,1]\to\frak H$ и $J\colon\frak H\to
W_{s,B,V}^{-m}[0,1]$. В этом слу\-чае оператор $T$ задаёт на пространстве $\frak H$
линейное отношение $T^\bullet\rightleftharpoons J^{-1}TI^{-1}$ с графиком
$$
	\{(y,z)\in\frak H\times\frak H\;:\;(\exists Y\in W_{s,A,U}^n[0,1])\quad
		(IY=y)\mathbin{\&}(TY=Jz)\}.
$$
Если при некотором $\lambda\in\Bbb C$ оператор $T-\lambda JI$ является ограниченно
обратимым, то очевидным образом определена также и ограниченная резольвента
$$
	(T^\bullet-\lambda)^{-1}=I\cdot(T-\lambda JI)^{-1}J.\leqno(\equation)
$$\label{eq:RezRil}%
Соответственно, график отношения $T^\bullet$ в этом случае заведомо замкнут.
Если дополнительно резольвента~\eqref{eq:RezRil} инъективна и имеет плотный образ,
то отношение $T^\bullet$ представляет собой оператор с плотной областью определения.

\subsubsection
{\it Если вложения $I\colon W_{s,A,U}^n[0,1]\to\frak H$ и $J\colon\frak H\to
W_{s,B,V}^{-m}[0,1]$ инъективны и име\-ют плотные образы, то допускающая
представление~\eqref{eq:RezRil} резольвента отношения $T^\bullet$ также инъективна
и имеет плотный образ.
}

\subsubsection\label{prop:nein}
{\it Пусть вложение $J\colon\frak H\to W_{s,B,V}^{-m}[0,1]$ инъективно, и пусть
для некоторого ограни\-ченного оператора $K\colon W_{1,A,U}^n[0,1]\to
C_{B,V}^m[0,1]$ в случае $s=1$ или {\relpenalty 5000 $K\colon W_{s,A,U}^n[0,1]\to
W_{s/(s-1),B,V}^m[0,1]$} в случае $s\neq 1$ справедливо тождество
$\langle z,IY\rangle\equiv\langle Jz,KY\rangle$. Пусть также равенство $\langle TY,
KY\rangle=0$ возможно лишь в случае $Y=0$. Тогда линейное отношение $T^\bullet$
представляет собой оператор. Если при этом индекс оператора $T$ равен нулю,
то существует ограниченный оператор $T^{-1}$, а область определения оператора
$T^\bullet$ плотна в $\frak H$.
}%
\proof
Если пара $(0,z)\in\frak H\times\frak H$ принадлежит графику отношения $T^\bullet$,
то должна существовать вектор-функция $Y\in W_{s,A,U}^n[0,1]$ со свойствами $TY=Jz$
и $IY=0$. При этом, однако, должны выполняться также равенства
$$
	\eqalign{\langle TY,KY\rangle&=\langle Jz,KY\rangle\cr
		&=\langle z,IY\rangle\cr &=0,}
$$
согласно сделанным предположениям гарантирующие справедливость равенства $Y=0$,
а потому и равенства $z=0$.

Далее, всякая вектор-функция $Y\in W_{s,A,U}^n[0,1]$ со свойством $TY=0$ заведомо
удовлетворяет равенству $\langle TY,KY\rangle=0$, а потому является нулевой.
Соответственно, факт ограни\-ченной обратимости оператора $T$ в случае равенства
его индекса нулю немедленно вытекает из альтернативы Фредгольма.

Наконец, если вектор $z\in\frak H$ ортогонален образу оператора $IT^{-1}J$,
то должны быть справедливы равенства
$$
	\eqalign{0&=\langle z,IT^{-1}Jz\rangle\cr
		&=\langle Jz,KT^{-1}Jz\rangle\cr
		&=\langle T\cdot[T^{-1}Jz],K\cdot[T^{-1}Jz]\rangle,}
$$
согласно сделанным предположениям гарантирующие выполнение равенства $T^{-1}Jz=0$,
а потому и равенства $z=0$.
\endproof

\subsubsection\label{prop:sekt}
{\it Пусть для некоторого ограниченного оператора $K\colon W_{1,A,U}^n[0,1]\to
C_{B,V}^m[0,1]$ в случае $s=1$ или $K\colon W_{s,A,U}^n[0,1]\to
W_{s/(s-1),B,V}^m[0,1]$ в случае $s\neq 1$ справедливо тождество
$\langle z,IY\rangle\equiv\langle Jz,KY\rangle$, и пусть числовая область значений
оператора $K^*T$ включается в некоторый сектор комплексной плоскости, имеющий
вершиной точку $0$. Тогда числовая область значений отношения $T^\bullet$
включается в тот же сектор.
}%
\proof
Если пара $(y,z)\in\frak H\times\frak H$ принадлежит графику отношения $T^\bullet$,
то должна существовать вектор-функция $Y\in W_{s,A,U}^n[0,1]$ со свойствами $IY=y$
и $TY=Jz$. Выполняющиеся при этом равенства
$$
	\eqalign{\langle z,y\rangle&=\langle Jz,KY\rangle\cr
		&=\langle TY,KY\rangle}
$$
как раз и гарантируют справедливость доказываемого утверждения.
\endproof

\subsubsection\label{prop:sim}
{\it Пусть вложение $I\colon W_{s,A,U}^n[0,1]\to\frak H$ инъективно, для некоторого
ограниченного оператора $K\colon W_{1,A,U}^n[0,1]\to C_{B,V}^m[0,1]$ в случае $s=1$
или {\relpenalty 5000 $K\colon W_{s,A,U}^n[0,1]\to W_{s/(s-1),B,V}^m[0,1]$}
в случае $s\neq 1$ справедливо тождество $\langle z,IY\rangle\equiv
\langle Jz,KY\rangle$, а оператор $K^*T$ является симметрическим. Тогда при любом
невещественном $\lambda\in\Bbb C$ оператор $T-\lambda JI$, если является
фредгольмовым с нулевым индексом, обладает ограниченным обратным.
}%
\proof
В рассматриваемом случае для всякой вектор-функции $Y\in W_{s,A,U}^n[0,1]$
справедливо равенство
$$
	\Im\langle (T-\lambda JI)Y,KY\rangle=-\Im\lambda\cdot\|IY\|^2.
$$
Соответственно, уравнение $(T-\lambda JI)Y=0$ не может иметь нетривиальных решений,
что, согласно альтернативе Фредгольма, как раз и означает справедливость
доказываемого утверждения.
\endproof

\medskip
Отметим в заключение, что, согласно утверждению~\S$\,$2.1.3
работы~[\cite{Vl:2014}], в случае инъек\-тивности и плотности вложений $I\colon
W_{s,A,U}^n[0,1]\to\frak H$ и $J\colon\frak H\to W_{s,B,V}^{-m}[0,1]$, а также
непустоты резольвентного множества пучка $T^\natural\colon\lambda\mapsto
T-\lambda JI$, спектр оператора $T^\bullet$ в точности совпадает со спектром
пучка $T^\natural$.


\section{Примеры}\label{par:Ekz}
\subsection
В качестве первого примера рассмотрим задачу Дирихле для дифференциального
уравнения
$$
	(py'')''-(q'y')'+r''y=f\in L_2[0,1],
$$
где $p,p^{-1}\in L_\infty[0,1]$ и $q,r\in L_2[0,1]$ (классический случай, когда
коэффициенты $p$, $q'$ и $r''$ предполагаются гладкими, связан с существенно более
сильными ограничениями). Этой задаче отвечает оператор $T\colon\Wo_2^2[0,1]\to
\Wo_2^{-2}[0,1]$ вида
$$
	\langle Ty,z\rangle\equiv\int_0^1\left[(py''-qy'+ry)\,\overline{z''}+
		(-qy''+2ry')\,\overline{z'}+ry''\overline{z}\right]\,dx.
$$
Данный оператор имеет форму~\eqref{eq:kvdif}, и потому уравнение $T^\bullet y=f$
для соответ\-ствующего неограниченного оператора в пространстве $L_2[0,1]$
равносильно граничной задаче
$$
	\belowdisplayskip 1\jot
	\eqalignno{(y^{[0]})'&=y^{[1]},\cr (y^{[1]})'&=-p^{-1}ry^{[0]}+
		p^{-1}qy^{[1]}+p^{-1}y^{[2]},\cr (y^{[2]})'&=p^{-1}qry^{[0]}+
		(2r-p^{-1}q^2)y^{[1]}-p^{-1}qy^{[2]}-y^{[3]},\cr
		(y^{[3]})'&=-p^{-1}r^2y^{[0]}+p^{-1}qry^{[1]}+p^{-1}ry^{[2]}-f,}
$$
$$
	\abovedisplayskip 0pt
	y^{[0]}(0)=y^{[1]}(0)=y^{[0]}(1)=y^{[1]}(1)=0.
$$
Указанная система уравнений, ожидаемым образом, в точности совпадает с системой
из леммы~2 работы~[\cite{Vl:2004}].

Развитая выше теория, однако, может быть приложена и к случаю коэффициентов
существенно более широкого вида. В частности, пусть коэффициент $p^{-1}$ есть
обобщённая производная некоторой строго возрастающей функции $H\in C[0,1]$,
коэффициент $q$ принадлежит классу $L_2([0,1];\,H')$, а коэффициент $r$~---
классу $L_2([0,1];\,H')\cap L_1[0,1]$. Тогда могут быть введены в рассмотрение
функции $\xi,\eta\in W_1^1[0,1]$ со свойствами
$$
	\xi\left({x+H(x)-H(0)\over 1+H(1)-H(0)}\right)\equiv x,\qquad
	\eta\left({x+H(x)-H(0)\over 1+H(1)-H(0)}\right)\equiv H(x),
$$
а также связанная с ними система суммируемых функций
$$
	A_{ij}\rightleftharpoons\cases{\xi'&при $i=0$, $j=1$,\cr
		\eta'&при $i=1$, $j=2$,\cr 0&иначе.}
$$
Рассматриваемая задача получит при этом трактовку $T^\bullet y=f$,
где действующий в про\-странстве $L_2[0,1]$ неограниченный оператор $T^\bullet=
(I^*)^{-1}TI^{-1}$ строится на основе инъективного вложения $I\colon
\Wo_{2,A}^2[0,1]\to L_2[0,1]$ вида
$$
	[IY](x)\equiv Y_0\left({x+H(x)-H(0)\over 1+H(1)-H(0)}\right)
$$
и оператора $T\colon\Wo_{2,A}^2[0,1]\to\Wo_{2,A}^{-2}[0,1]$ вида
$$
	\langle TY,Z\rangle\equiv\int_0^1\left[\eta'\cdot(Y_2-\sigma Y_1+
		\rho Y_0)\,\overline{Z_2}+(-\eta'\sigma Y_2+2\xi'\rho Y_1)\,
		\overline{Z_1}+\eta'\rho Y_2\,\overline{Z_0}\right]\,dx,
$$
где положено $\sigma\rightleftharpoons q\circ\xi$ и $\rho\rightleftharpoons
r\circ\xi$. Результаты предыдущего параграфа позволяют теперь представить исходную
задачу в форме
$$
	\belowdisplayskip 1\jot
	\eqalignno{(y^{[0]})'&=\xi'y^{[1]},\cr (y^{[1]})'&=-\eta'\rho y^{[0]}+
		\eta'\sigma y^{[1]}+\eta'y^{[2]},\cr
		(y^{[2]})'&=\eta'\sigma\rho y^{[0]}+(2\xi'\rho-
		\eta'\sigma^2)y^{[1]}-\eta'\sigma y^{[2]}-\xi'y^{[3]},\cr
		(y^{[3]})'&=-\eta'\rho^2y^{[0]}+\eta'\sigma\rho y^{[1]}+
		\eta'\rho y^{[2]}-f\circ\xi,}
$$
$$
	\abovedisplayskip 0pt
	y^{[0]}(0)=y^{[1]}(0)=y^{[0]}(1)=y^{[1]}(1)=0,
$$
где $y^{[0]}=y\circ\xi$.

\subsection
В качестве второго примера возьмём задачу периодического типа
для дифферен\-циального уравнения третьего порядка
$$
	-iy'''-(py')'+q'y=f\in L_2[0,1],\leqno(\equation)
$$\label{eq:tria_ordo}%
где $p,q\in L_1[0,1]$. В случае задач нечётного порядка естественным является
рассмо\-трение ситуации $s=1$, в отличие от типичного для чётного случая равенства
$s=2$. Соответственно, свяжем с задачей пару пространств
$$
	\eqalignno{W_{1,A,U}^2&\rightleftharpoons\{y\in W_1^2[0,1]\;:\;
			y(1)-y(0)=y'(1)-y'(0)=0\},\cr
		C_{B,V}^1&\rightleftharpoons\{y\in C^1[0,1]\;:\;y(1)-y(0)=0\}.
	}
$$
Умножая обе стороны уравнения~\eqref{eq:tria_ordo} на компоненту $Z_0$ произвольной
вектор-функции $Z\in C_{B,V}^1[0,1]$, интегрированием по частям убеждаемся
в справедливости равенства $T^\bullet y=f$, где отношение $T^\bullet$ отвечает
оператору $T\colon W_{1,A,U}^2[0,1]\to W_{1,B,V}^{-1}[0,1]$ вида
$$
	\langle TY,Z\rangle\equiv\int_0^1\left[(iY_2+pY_1-qY_0)\,\overline{Z_1}-
		qY_1\,\overline{Z_0}\right]\,dx.
$$
Результаты предыдущего параграфа дают теперь для уравнения $T^\bullet y=f$
эквивалентное представление
$$
	\belowdisplayskip 1\jot
	\eqalignno{(y^{[0]})'&=y^{[1]},\cr (y^{[1]})'&=-iqy^{[0]}+ipy^{[1]}
		-iy^{[2]},\cr (y^{[2]})'&=-qy^{[1]}-f,}
$$
$$
	\abovedisplayskip 0pt
	y^{[0]}(1)-y^{[0]}(0)=y^{[1]}(1)-y^{[1]}(0)=y^{[2]}(1)-y^{[2]}(0)=0.
$$
На ос\-нове утверждений~\ref{prop:sekt} и~\ref{prop:sim}, с учётом наличия
ограниченного вложения прост\-ранства $W_{1,A,U}^2[0,1]$ в пространство
$C_{B,V}^1[0,1]$, легко показывается также, что в случае вещественнозначности
коэффициентов $p,q\in L_1[0,1]$ неограниченный оператор $T^\bullet$ является
самосопряжённым.

\subsection
В качестве последнего примера рассмотрим задачу Дирихле для исследовавшегося
в рабо\-тах~[\cite{Fr:2011}] и~[\cite{Vl:2012}] уравнения
$$
	-{d\over dG}\left({dy\over dH}\right)=f\in L_2([0,1];\,G'),
$$
где $H\in C[0,1]$ есть некоторая неубывающая функция со свойствами $H(0)=H(1)-1=0$,
а неубывающая функция $G\in C[0,1]$ допускает представление $G(x)\equiv N(H(x))$.
В этом случае решение задачи заведомо  может быть записано в виде $y=Iu$,
где вложение $I\colon\Wo_2^1[0,1]\to L_2([0,1];\,G')$~--- вообще говоря,
не инъективное~--- подчиняется тождеству $[Iu](x)\equiv u(H(x))$. Кроме того,
заведомо найдётся функция $g\in L_2([0,1];\,N')$, для которой композиция $g\circ H$
будет совпадать в пространстве $L_2([0,1];\,G')$ с исходной правой частью $f$,
а потому и подчиняться тождеству
$$
	(\forall v\in\Wo_2^1[0,1])\qquad\langle I^*f,v\rangle=
		\int_0^1 g\overline{v}\,dN.
$$
Соответственно, исходная задача допускает переформулировку $(I^*)^{-1}TI^{-1}y=f$,
где $T\colon\Wo_2^1[0,1]\to\Wo_2^{-1}[0,1]$ есть обычный оператор двукратного
дифференцирования
$$
	\langle Tu,v\rangle\equiv\int_0^1 u'\overline{v'}\,dx.
$$
Утверждения~\ref{prop:nein} и~\ref{prop:sekt} вместе с очевидным фактом
положительности оператора $T$ позволяют гарантировать корректную определённость,
самосопряжённость и положи\-тельность действующего в пространстве $L_2([0,1];\,G')$
неограниченного оператора $T^\bullet\rightleftharpoons (I^*)^{-1}TI^{-1}$.


\vskip 0.4cm
\eightpoint\rm
{\leftskip 0cm\rightskip 0cm plus 1fill\parindent 0cm
\bf Литература\par\penalty 20000}\vskip 0.4cm\penalty 20000
\bibitem{Vl:2004} А.$\,$А.~Владимиров. {\it О сходимости последовательностей
обыкновенных дифференциальных опе\-раторов}~// Матем. заметки.~--- 2004.~---
Т.~75, \No~6.~--- С.~941--943.

\bibitem{Vl:2016} А.$\,$А.~Владимиров. {\it К вопросу об осцилляционных свойствах
положительных дифференциаль\-ных операторов с сингулярными коэффициентами}~//
Матем. заметки.~--- 2016.~--- Т.~100, \No~6.~--- С.~800--806.

\bibitem{NSh:1999} М.$\,$И.~Нейман--заде, А.$\,$А.~Шкаликов. {\it Операторы
Шрёдингера с сингулярными потенциалами из прост\-ранств мультипликаторов}~//
Матем. заметки.~--- 1999.~--- Т.~66, \No~5.~--- С.~723--733.

\bibitem{MSh:2016} К.$\,$А.~Мирзоев, А.$\,$А.~Шкаликов. {\it Дифференциальные
операторы чётного порядка с коэффициен\-тами-распределениями}~//
Матем. заметки.~--- 2016.~--- Т.~99, \No~5.~--- С.~788--793.

\bibitem{Vl:2014} А.$\,$А.~Владимиров. {\it Теоремы о представлении и вариационные
принципы для самосопряжённых опера\-торных матриц}~// arXiv:1403.2253.

\bibitem{Fr:2011} U.~Freiberg. {\it Refinement of the spectral asymptotics
of generalized Krein Feller operators}~// Forum Math.~--- 2011.~--- V.~23.~---
P.~427--445.

\bibitem{Vl:2012} А.$\,$А.~Владимиров. {\it Об одном классе сингулярных задач
Штурма--Лиувилля}~// arXiv:1211.2009.
\bye